\documentclass[number,citesort,dvips]{arxbj}
\usepackage{upgreek}
\usepackage{graphicx}

\aid{0}
\volume{17}
\issue{4}
\pubyear{2011}
\firstpage{1420}
\lastpage{1434}
\doi{10.3150/10-BEJ322}

\makeatletter
\newtheorem{theorem}{Theorem}
\newtheorem{lemma}{Lemma}

\newcommand{\mvec}{\mathbf{m}}
\newcommand{\xvec}{\mathbf{x}}

\newcommand{\Xvec}{\mathbf{X}}

\newcommand{\hvec}{\mathbf{h}}

\newcommand{\zvec}{\mathbf{z}}
\newcommand{\yvec}{\mathbf{y}}

\newcommand{\muvec}{\bolds{\mu}}
\newcommand{\alphavec}{\bolds{\alpha}}

\newcommand{\sigmavec}{\bolds{\sigma}}
\makeatother

\begin{document}
\begin{frontmatter}

\title{Some intriguing properties of Tukey's half-space depth}
\runtitle{Tukey's half-space depth}

\begin{aug}
\author{\fnms{Subhajit}  \snm{Dutta}\thanksref{e1}\ead[label=e1,mark]{subhajit\_r@isical.ac.in}},
\author{\fnms{Anil K.} \snm{Ghosh}\thanksref{e2}\ead[label=e2,mark]{akghosh@isical.ac.in}}
  \and
\author{\fnms{Probal} \snm{Chaudhuri}\thanksref{e3}\corref{}\ead[label=e3,mark]{probal@isical.ac.in}}
\runauthor{S. Dutta, A.K. Ghosh and P. Chaudhuri}
\address{Theoretical Statistics and Mathematics Unit, Indian Statistical Institute, Kolkata 700108,
 India.\break \printead{e1,e2,e3}}
\end{aug}

\received{\smonth{6} \syear{2010}}

\begin{abstract}
For multivariate data, Tukey's half-space depth is one of the most
popular depth functions available in the literature. It is conceptually
simple and satisfies several desirable properties of depth functions.
The Tukey median, the multivariate median associated with the
half-space depth, is also a well-known measure of center for
multivariate data with several interesting properties. In this article,
we derive and investigate some interesting properties of half-space
depth and its associated multivariate median. These properties, some of
which are counterintuitive, have important statistical consequences in
multivariate analysis. We also investigate a natural extension of
Tukey's half-space depth and the related median for probability
distributions on any Banach space (which may be finite- or
infinite-dimensional) and prove some results that demonstrate anomalous
behavior of half-space depth in infinite-dimensional spaces.
\end{abstract}


\end{frontmatter}

\section{Introduction}\label{sec1}

Over the last three decades, data depth has emerged as a powerful
concept leading to the generalization of many univariate statistical
methods to the multivariate setup. A depth function measures the
centrality of a point $\xvec$ with respect to a data set or a
probability distribution and thus helps to define an ordering and a
version of ranks for multivariate data. There are several notions of
data depth available in the literature (see, e.g.,
\cite{Liu1999,Zuo2000,Vardi,Mosler,Mizera,Lopez}). Tukey's half-space
depth (see~\cite{Tukey}) is one of the most popular depth functions
used by many researchers. The construction of central regions based on
trimming (see, e.g., \cite{Nolan}), robust estimation of multivariate
location (see, e.g., \cite{Donoho}), tests of multivariate statistical
hypotheses (see, e.g., \cite{Chaudhuri}) and supervised classification
(see, e.g., \cite{Ghosh2005a,Ghosh2005a}) are some examples of its
widespread application.

Like other popular depth functions, half-space depth has some nice
theoretical properties. In fact, it satisfies all four of the desirable
properties of depth functions first mentioned in~\cite{Liu1990} and
subsequently investigated in \cite{Zuo2000}, namely, affine invariance,
maximality at the center, monotonicity with respect to the deepest
point and vanishing at infinity. Moreover, if the underlying population
distribution $F$ has a spherically symmetric density $f$, that is,
$f(\xvec)=\psi(\|\xvec\|_2)$ for some $\psi \dvtx \mathbb{R}_+
\rightarrow \mathbb{R}_+$, the half-space depth turns out to be a
decreasing function of $\|\xvec\|_2= (|x_1|^2+\cdots+|x_d|^2)^{1/2}$.
Consequently, when $\psi$ is monotonically decreasing (i.e., $f$ is
unimodal), the half-space depth becomes an increasing function of $f$
and vice versa.
Therefore, in such cases, the half-space depth contours coincide with
the contours of the density function. Because of this property of
the
half-space depth, classification rules based on the ordering of the
half-space depth functions coincide with the optimal\vadjust{\goodbreak} Bayes classifier
for discriminating among spherically symmetric unimodal populations
differing in their centers of symmetry (see, e.g., \cite{Ghosh2005b}).
Similarly, the use of the half-space depth functions to order and trim
multivariate data sets (see, e.g., \cite{Nolan,Donoho}) leading to the
determination of central and outlying observations has a natural
justification when the density contours coincide with the half-space
depth contours. Also, due to this relation between half-space depth and
spherical symmetry, half-space depth has been used to construct
diagnostic tools for checking spherical symmetry of a data cloud (see,
e.g., \cite{Liu1999}, pages 809--811). Another well-known feature of
half-space depth is its characterization property. Koshevoy \cite{Koshevoy2002}
proved that if the half-space depth functions of two atomic measures
with finite support are identical, then the measures are also
identical. Cuesta-Albertosa and  Nieto-Reyes \cite{Cuesta} proved this characterization property of Tukey
depth for discrete distributions. Under some regularity conditions,
Koshevoy \cite{Koshevoy2003} proved this characterization property for
absolutely continuous probability distributions with compact support in
finite-dimensional spaces. Hassairi and Regaieg~\cite{Hassairi} generalized it to absolutely
continuous distributions with connected supports.

However, the half-space depth function  has several limitations.
The half-space median derived from half-space depth has a lower
breakdown point and relative efficiency compared to the median based on
projection depth (see \cite{Zuo2003}). Dang and Serfling \cite{Dang} pointed out that the
outlier identifier based on the half-space depth has a ``severe'' and
``unacceptable'' trade-off between ``masking breakdown point'' and
``false positive rate''. Moreover, if the half-space depth contours fail
to match the density contours, then the classifiers based on half-space
depth may lead to misclassification rates higher than the Bayes risk.
The diagnostic tool developed in \cite{Liu1999}, pages 809--811 for
detecting deviations from spherical symmetry using half-space depth
also relies heavily on the fact that under $l_2$-symmetry, the depth
contours are concentric spheres with half-space median at the center.
So, in the absence of this property of the half-space depth contours,
such a diagnostic tool may not lead to useful results. Now, a natural
question that arises from this discussion is whether this property of
half-space depth contours holds for other symmetric distributions, for
example, in the case of $l_p$-symmetric distributions, when
$f(\xvec)=\psi(\|\xvec\|_p)$ for some $p \neq 2$ and $\psi$ is
monotonically decreasing. Here, for any $p>0$ and
$\xvec=(x_1,\ldots,x_d) $ $\in \mathbb{R}^d$, we define $\|\xvec\|_p =
(|x_1|^p+\cdots+|x_d|^p)^{1/p}$. In Section \ref{sec2}, we carry out an
investigation to answer this question.

For any continuous univariate distribution, it is straightforward to
see that the median is the point with  half-space depth 0.5. In Section
\ref{sec3}, we investigate to what extent this property of half-space median
holds for multivariate continuous distributions and derive a
characterization of the multivariate distribution for which the
half-space depth of Tukey median will achieve its maximum value, namely
0.5. We propose a statistical test for angular symmetry of continuous
multivariate distributions based on this characterization and briefly
study the performance of the proposed test. In this section, we also
consider natural extensions of half-space depth and half-space median
for probability distributions in arbitrary Banach spaces using the
concept of linear functionals on such spaces. Some anomalous behaviors
of half-space depth for probability distributions on
infinite-dimensional spaces and their implications are discussed in
Section \ref{sec4}. Proofs of theorems and lemmas (along with their statements)
are deferred to the \hyperref[appendix]{Appendix}.

\section{Half-space depth contours for $l_p$-symmetric density
functions}\label{sec2}

In this section, we study the behavior of the half-space depth contours
for a wide class of symmetric distributions. As was mentioned in the
\hyperref[sec1]{Introduction}, the half-space depth contours coincide with the density
contours if the p.d.f. $f$ is such that $f(\xvec)=\psi(\|\xvec\|_2)$
for some monotonically decreasing $\psi \dvtx \mathbb{R}_{+} \rightarrow
\mathbb{R}_{+}$, and this is an important feature of half-space depth
with many useful statistical applications. Here, we will investigate
the situation when $\| \cdot \|_2$ is replaced by $\| \cdot \|_p$,
where $p$ is positive and $p \neq 2$.

\subsection{Depth contours for $p=\infty$}

\begin{figure}

\includegraphics{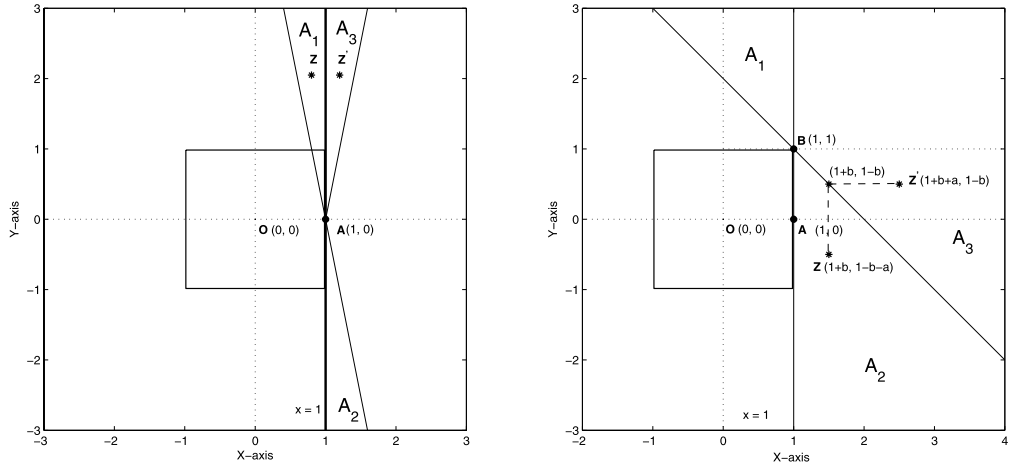}

  \caption{$l_\infty$ contour and the line defining
the half-space depth of ($1,0$) and ($1,1$).}\label{fig1}
\end{figure}

For $p=\infty$, the
p.d.f. $f(\xvec)=f(x_1,x_2,\ldots,x_d)=\psi(\max\{|x_1|,|x_2|,\ldots,|x_d|\})$
for some monotonically decreasing function $\psi$. Clearly, the density
contours here are concentric $d$-dimensional hypercubes with the origin at
the center. We now check whether or not all points on the surface of a
hypercube with origin at the center have the same depth. First,
consider the point $A=(1,0,\ldots,0)$ on the surface of the unit
hypercube $\{ \xvec \dvtx \|\xvec\|_{\infty} =1\}$ (see Figure \ref{fig1} for a~diagram in the case $d=2$). It can be shown that the hyperplane $x_1=1$
determines the half-space depth of this point, and this depth is
$P(X_1\ge 1)$, where $\Xvec=(X_1,X_2,\ldots,X_d)$ has the
p.d.f. $f(\xvec)$ (see Lemma \ref{lem1} in the \hyperref[appendix]{Appendix}).

Note that the line $x_1=1$ also passes through the point
$B=(1,1,0,\ldots,0)$ (see the right-hand diagram in Figure \ref{fig1} when
$d=2$). So, $A$ and $B$ will have the same depth if and only if there
exists no other hyperplane that passes though $B$ in such a way that
the probability of one of its half-spaces is smaller than $P(X_1 \ge
1)$. However, the hyperplane $x_1+x_2=2$ passes through the point~$B$,
and we can show that $P(X_1+X_2 \ge 2)< P(X_1 \ge 1)$ (see Lemma \ref{lem2} in
the \hyperref[appendix]{Appendix}). This implies that if the p.d.f. $f$ is of the form
$f(\xvec) = \psi(\|\xvec\|_\infty)$ with a monotonically decreasing~$\psi$, then the half-space depth contours cannot coincide with the
corresponding density contours.

\begin{figure}

\includegraphics{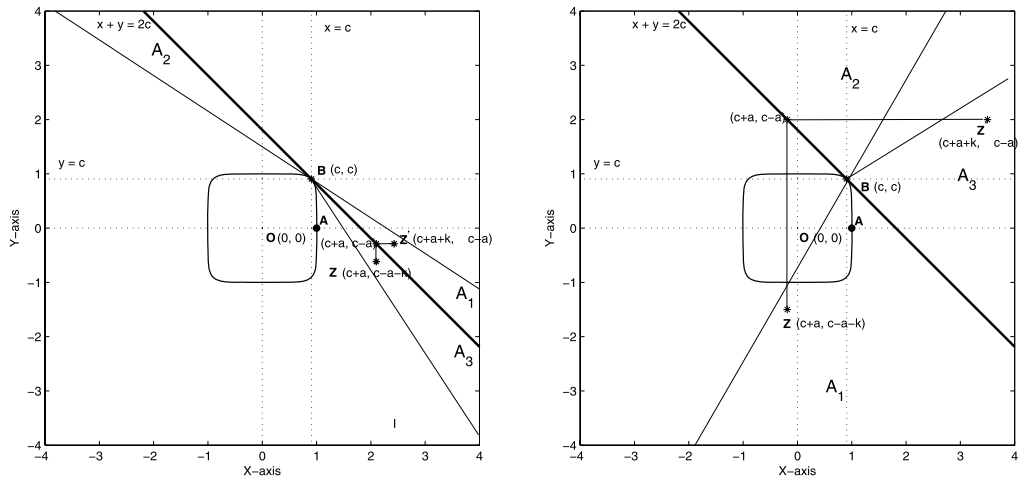}

  \caption{$l_p$ contour and the lines defining the
half-space depth of ($c,c$) for $p=5$.}\label{fig2}
\end{figure}

\subsection{Depth contours for $1 \le p <\infty$}

Next, consider the case where $1\le p<\infty$. Clearly,
$A=(2^{1/p}c,0,0,\ldots,0)$ and $B=(c,c,\break 0,\ldots,0)$ are two points on
the same $l_p$ contour (see Figure \ref{fig2} for the case $d=2$). First, we
check whether or not the half-space depths of these two points are
equal. In view of Lemma \ref{lem1}, the depth of $A$ is given by $P(X_1\ge
2^{1/p}c)$ when $c>0$.
We can also prove that the hyperplane $x_1+x_2=2c$ determines the
half-space depth of $B$ and that this depth is $P(X_1+X_2 \ge 2c)$ (see
Lemma \ref{lem3} in the \hyperref[appendix]{Appendix}).

It follows from the discussion in the preceding paragraph that the two
points $A$ and $B$ will have the same depth only if $P(X_1 \ge
2^{1/p}c)=P(X_1+X_2 \geq 2c)$. Note that here  we can choose $c$
arbitrarily. Therefore, the depth and the density contours can coincide
only if $P(X_1 \ge 2^{1/p}c)=P(X_1+X_2 \geq 2c)$ for all values of $c$,
that is, only if $X_1$ and $2^{\alpha}(X_1+X_2)$ are identically
distributed for $\alpha=(1-p)/p$. Now, if we assume the existence of
the second order moments of the $X_i$'s, then the equality of the
variances of $X_1$ and $2^{\alpha}(X_1+X_2)$ and the fact that $X_1$
and $X_2$ are uncorrelated (in view of the $l_p$-symmetry of the
density $f$) imply that $\alpha=-1/2$ or $p=2$. Even if we do not
assume any moment condition, the above result holds (see Lemma \ref{lem4} in the
\hyperref[appendix]{Appendix}). Also, it is interesting to note that for $p<2$, we can
always choose a $c$ such that the depth of $B$ is more than that of
$A$. On the other hand, for $p>2$, it is always possible to choose a~$c$ such that $A$ has larger depth than
$B$.\looseness=-1

\subsection{Depth contours for $p<1$}

\begin{figure}[b]

\includegraphics{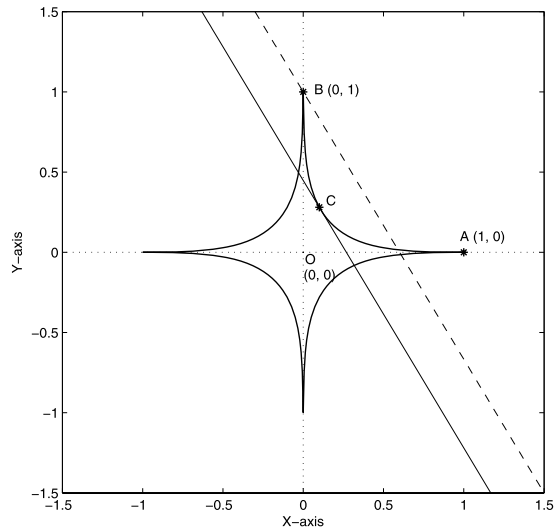}

  \caption{$l_p$ contour for the case   $p = 1/2$.}\label{fig3}
\end{figure}

Finally, we investigate the case $p<1$. Note that in this case, the
regions bounded by $l_p$ contours are no longer convex sets (see Figure
\ref{fig3} for the case $d=2$). Consider three points $A=(1,0,\ldots,0)$,
$B=(0,1,0,\ldots,0)$ and $C=(\alpha, \beta, 0, \ldots,0)$ on the same
$l_p$ contour, where $\alpha, \beta >0$ and $| \alpha |^p + | \beta |^p
= 1$. Consider any hyperplane passing through $C$. It will split
$\mathbb{R}^d$ into two half-spaces, one of which will contain the
origin. Since $p<1$, at least one of the two points $A$ and $B$ will
lie in the half-space that does not contain the origin. Without loss of
generality, we can assume that the hyperplane that determines the
half-space depth of $C$ puts $B$ and the origin in two different
half-spaces (see the bold line in Figure \ref{fig3} for the case $d=2$). We can
now make a~parallel shift of that hyperplane away from the origin until
it hits the point $B$ (see the dotted line in Figure \ref{fig3} for the case
$d=2$). Clearly, the half-space created by this new hyperplane that has
smaller probability measure will have smaller probability than that  of
each of the two half-spaces created by the older hyperplane. Therefore,
the half-space depth of $B$ has to be smaller than that of $C$ and
hence the depth contours cannot coincide with the density contours.

\begin{figure}

\includegraphics{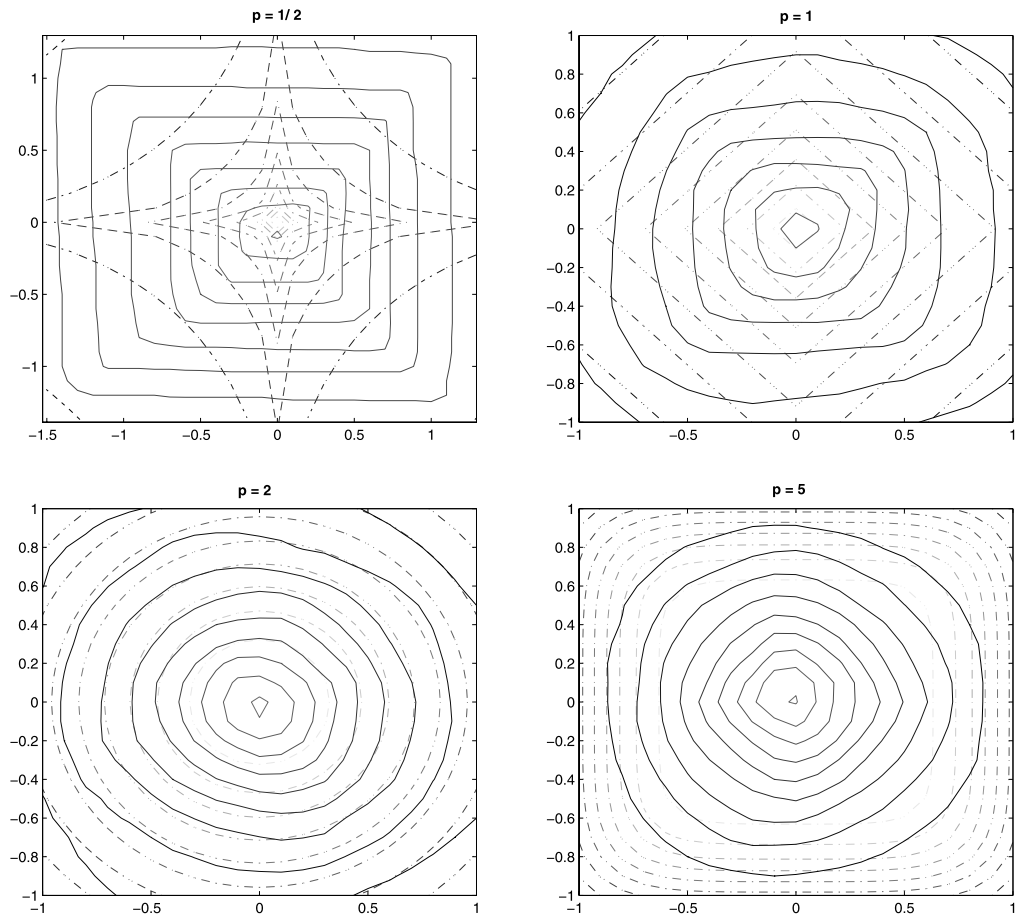}

  \caption{Density contours and their corresponding
half-space depth contours.}\label{fig4}
\vspace*{-9pt}
\end{figure}

Summarizing our discussion in this section, we now have the following
theorem.

\begin{theorem}\label{theo1}
Consider a probability
distribution on $\mathbb{R}^d$ with the p.d.f. $f$ such that
$f(\xvec)=\psi(\|\xvec\|_p)$ for some monotonically decreasing function
$\psi$. The half-space depth contours associated with $f$ will then
coincide with the density contours if and only if $p=2$.
\end{theorem}

Figure \ref{fig4} presents the empirical half-space depth contours (indicated
using connected lines) computed using 500 observations from bivariate
$l_p$-symmetric distributions with different values of $p$ (i.e.,
$p=1/2,1,2,5$). In each case, we consider the density to be of the form
$f(\xvec) = \frac{(2\Gamma(1/p))^2}{p^2}\exp(-\{|x_1|^p+|x_2|^p\})$ and
the corresponding density contours are also plotted (indicated using
dotted lines) in Figure \ref{fig4}. From this figure, it is quite evident that
the half-space depth contours and the density contours are markedly
different when $p\neq 2$. So, unlike what was done by \cite{Liu1999}, pages 809--811, we cannot develop a diagnostic tool for
checking $l_p$-symmetry using half-space depth when $p\neq 2$.

\begin{figure}

\includegraphics{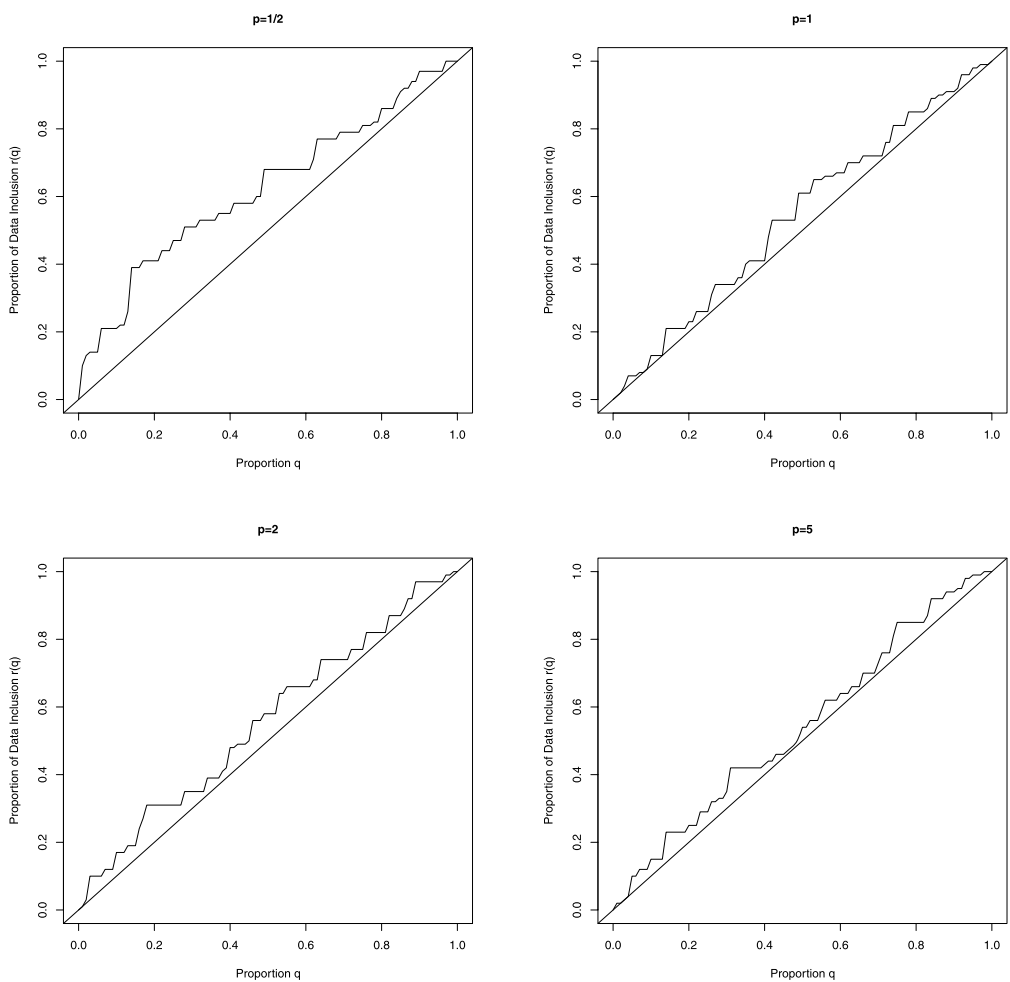}

  \caption{Diagnostic tool for checking spherical
symmetry.}\label{fig5}
\end{figure}

It is also of interest to note that along with $p=2$, for $p=1$ and
$5$, the half-space depth contours are nearly circular. Since the
diagnostic tool for spherical symmetry proposed in \cite{Liu1999},
pages~809--811, relies heavily on the sphericity of the depth contours,
it may fail to detect the deviation from spherical symmetry in the
cases $p=1$ and $5$. But for $p=1/2$, since the depth contours
are far from being circular, we can expect to detect this deviation
using their diagnostic tool. This is what we observed when we performed
the following experiment. Following \cite{Liu1999}, pages 809--811,
for different values of $q$ ($0<q<1$), we found the smallest sphere
$S_q$ containing the $q$th central hull and computed the fraction of
the data $r(q)$ lying in $S_q$. This fraction $r(q)$ is plotted against
$q$ for four different $l_p$-symmetric distributions with $p=1/2,1,2$
and $5$, and these plots are presented in Figure \ref{fig5}. Note that if the
underlying distribution is spherically symmetric (i.e.,
$l_2$-symmetric), the resulting curve should lie near the diagonal line
joining the points ($0, 0$) and ($1, 1$). The area between the curve
and the diagonal line gives an indication of the deviation from
spherical symmetry. As expected, for $p=1,2$ and $5$, these curves were
close to the diagonal line, but in the case $p=1/2$, the curve had a
significant deviation from the diagonal line (see Figure \ref{fig5}). So, the
diagnostic tool could detect the deviation from spherical symmetry only
in the case of $l_{1/2}$-symmetry.

We have seen that the half-space depth contours do not match  the
density contours for any $l_p$-symmetric distribution with $p\neq 2$,
and this leads to several limitations on statistical tools based on
half-space depth, as was already discussed in the \hyperref[sec1]{Introduction} and the
present section. However, it will be appropriate to note here that in
such cases, the depth function may provide some useful information
which may not be contained in the density function. While density is
only a local measure, which measures the local probability mass, depth
is a global measure, which gives useful information about global
features like the central and outlying points of a data cloud or
probability distribution. For instance, in the case of multivariate
uniform distributions, the density function, being constant, fails to
give any idea about the central and the peripheral points of the
distribution; however, the half-space depth function provides a
meaningful measure of central tendency, for example, by identifying the
point with the maximum depth (see~\cite{Serfling}).

\section{Half-space median and its depth}\label{sec3}

As we have already pointed out in the \hyperref[sec1]{Introduction}, for continuous
univariate distributions, the median is the point with half-space depth
0.5. In a sense, this is a very desirable and natural property for a
measure of the center of a distribution, and we would also like this
property to hold in a multivariate setup. If this property holds for a
multivariate distribution, any hyperplane passing through the median
will lead to two half-spaces having equal probability measures.
Unfortunately, as we will gradually see in this section, this may not
always be true for multivariate distributions, even if the distribution
is absolutely continuous with respect to the Lebesgue measure on a
Euclidean space.

Note that for any $l_p$-symmetric density function
$f(\xvec)=\psi(\|\xvec\|_p)$ with $0<p\le \infty$, the origin turns out
to be the half-space median with the half-space depth 0.5. In fact,
this is true whenever $\Xvec$ and $-\Xvec$ have the same distribution
(i.e., the distribution is centrally symmetric), or even under a
slightly weaker condition that any real-valued linear projection has
median zero. We should also note that in all these cases, the
half-space median coincides with the coordinatewise median, and the
depth of the half-space median, namely the origin, is 0.5. However,
this only holds for a special class of multivariate distributions. For
instance, for a bivariate uniform distribution on a right-angled
isosceles triangle, we can easily show that the half-space depth of any
point is smaller than 0.5. We can consider another interesting example
of a continuous bivariate distribution, where the p.d.f. $f$ has
support on $\{(x_1,x_2) \dvtx x_1+x_2\ge 0,  x_1x_2 \le 0\}$. In this case,
if $f$ is symmetric about the $x_1=x_2$ line, we can easily verify that
the half-space median will have depth smaller than 0.5, and the
coordinatewise median will have  zero half-space depth. We have already
indicated some sufficient conditions for the depth of the half-space
median to be 0.5, and in view of the two preceding examples, we would
like to know some necessary  {and} sufficient conditions for this.
We now state a theorem, the proof of which is given in the \hyperref[appendix]{Appendix}.

\begin{theorem}\label{theo2}
Suppose that $\Xvec$ is a
$d$-dimensional random  vector with a probability distribution which
has its half-space median at $\muvec \in \mathbb{R}^d$. Then, the
half-space depth of $\muvec$ will be 0.5 if and only if
$(\Xvec-\muvec)/\|\Xvec-\muvec\|_2$ and
$(\muvec-\Xvec)/\|\Xvec-\muvec\|_2$ are identically distributed.
\end{theorem}

This theorem implies that the half-space median will have depth 0.5 if
and only if the underlying distribution is angularly symmetric.
Liu et al. \cite{Liu1999}, pages 811--814, stated the sufficient part of this result
and used it to develop a diagnostic tool for verification of angular
symmetry of a distribution. This necessary and sufficient condition can
also be used to develop a statistical test for the angular symmetry of
a distribution. As discussed in \cite{Small},  Ajne's test (see \cite{Ajne}),
which is a~distribution-free test for bivariate data, can be
used for testing angular symmetry of a bivariate distribution about a
specified point (say, $\muvec_0$). However, the test that we propose
here is applicable to multivariate data in any dimension and does not
require any specification of the center of symmetry, which is estimated
from the data. Given a random sample $\xvec_1, \xvec_2, \ldots,
\xvec_n$ of size $n$, let $\tilde \mvec_n$ be the half-space median and
$\Delta_n$ denote the half-space depth of $\tilde \mvec_n$ in that
sample. For testing the null hypothesis of angular symmetry, an ideal
procedure would be to reject the null hypothesis if $\Delta_n < c_n$,
where $c_n$ is an appropriate percentile (that depends on the specified
level of the test) of the distribution of $\Delta_n$ under the null
hypothesis. However, it is not possible to determine an exact value of
$c_n$ in practice because the distribution of $\Delta_n$ depends on the
underlying angularly symmetric distribution of the data, which is
usually not specified in practice.
In practice, we propose that for a random sample $\xvec_1, \xvec_2,
\ldots, \xvec_n$, we first compute $\yvec_i=\xvec_i-\tilde \mvec_n$ for
$i=1,2,\ldots,n$, generate i.i.d. observations $z_1, z_2, \ldots , z_n$
such that $P(z_i = 1) = P(z_i = -1) = 1/2$ and then compute $\xvec_i^*
= z_i \yvec_i + \tilde \mvec_n$ for $i=1,2,\ldots,n$. This procedure is
motivated by the well-known idea of bootstrapping. These $\xvec_i^*$'s
can be viewed like a ``bootstrap sample'' generated from the original
sample under the null hypothesis of symmetry, and we can calculate the
depth $\Delta_n^*$ of the half-space median $\tilde \mvec_n^*$ based on
that ``bootstrap sample''. We can repeat this ``bootstrap procedure''
$M$ times depending on our computing resources and denote by $
\Delta_{n,m}^*$ the half-space depth of the half-space median in the
$m$th ``bootstrap sample'' ($m=1,2,\ldots,M$). The critical value $c_n$
mentioned earlier can then be estimated from the ``bootstrap empirical
distribution'' of $ \Delta_{n}^*$. In other words, for a specified
level $0 < \alpha < 1$, the null hypothesis of angular symmetry is to
be rejected if $\sum_{m=1}^{M} I\{ \Delta_{m,n}^* \le \Delta_n\}/M <
\alpha$.

To evaluate the performance of our proposed test, we carried out a
thorough simulation study with six examples using the software package
R. In each case, we generated samples of size 50 and 100, implemented
our test using $M=1000$ ``bootstrap samples'' and, in order to estimate
     the probability of rejection of $H_0$ by the test,  repeatedly applied it on 1000
Monte Carlo replications in dimensions $d={}$ 2, 3 and 4. The first five
examples were motivated by five bivariate examples in \cite{Liu1999}, page 814, which include three examples with angularly
symmetric distributions, namely
 D1, D2 and D3, and two examples, namely D4 and D5, where the underlying
 distributions were not angularly symmetric (\cite{Liu1999}, page 814, for a
 detailed description of these examples). Here, we consider the natural
 multivariate version of these five examples. In the last example, D6,
 which is also not angularly symmetric, when $d=2$, we generated observations
 from a bivariate uniform distribution on the right-angled isosceles triangle
 formed by the points $(0,0)$, $(1,0)$ and $(0,1)$. For an extension of D6 in
 dimensions $d > 2$, we have considered the simplex formed by the origin,
 the coordinate axes and the hyperplane $x_1 + \cdots + x_d = 1$ in $\mathbb{R}^d$ in
 place of the triangle. Table \ref{tab1} reports the proportion of cases, out of 1000 Monte
 Carlo replications, where the null hypothesis was rejected for two nominal
 values of $\alpha$, namely, 0.05 and 0.01. This table clearly shows good
 level as well as power properties of the proposed test procedure.

\begin{table}
\tabcolsep=0pt
\caption{Probability of rejection of $H_0$ by the proposed test}\label{tab1}
\begin{tabular*}{\textwidth}{@{\extracolsep{\fill}}llcccccccccccc@{}}
\hline
$d \downarrow$ & Data sets $\rightarrow$ & \multicolumn{2}{c}{D1} & \multicolumn{2}{c}{D2} &\multicolumn{2}{c}{D3} & \multicolumn{2}{c}{D4}
& \multicolumn{2}{c}{D5}  & \multicolumn{2}{c@{}}{D6} \\ [-7pt]
& & \multicolumn{2}{c}{\hrulefill} & \multicolumn{2}{c}{\hrulefill} &\multicolumn{2}{c}{\hrulefill} & \multicolumn{2}{c}{\hrulefill}
& \multicolumn{2}{c}{\hrulefill}  & \multicolumn{2}{c@{}}{\hrulefill}\\
[-3pt]
& Nominal $\rightarrow$ &$1\%$ &$5\%$ &$1\%$ &$5\%$ &$1\%$ &$5\%$ &$1\%$ &$5\%$
&$1\%$ &$5\%$ &$1\%$ &$5\%$ \\
& level ($\alpha$) & & & & & & & & & & & & \\
\hline
2 &$n=50$ &0.012 &0.052 &0.012 &0.054 &0.010 &0.044 &0.170 &0.318 &0.406 &0.663 &0.247 &0.418 \\
 &$n=100$ &0.014 &0.054 &0.014 &0.053 &0.010 &0.058
&0.486 &0.728 &0.870 &0.960 &0.641 &0.846 \\ [5pt]
3 &$n=50$ &0.011
&0.044 &0.003 &0.035 &0.015 &0.057 &0.294 &0.554 &0.751 &0.869 &0.403
&0.662 \\
 &$n=100$ &0.009 &0.051 &0.006 &0.040 &0.012
&0.046 &0.822 &0.949 &0.996 &1.000 &0.929 &0.982 \\ [5pt]
 4 &$n=50$
&0.009 &0.054 &0.013 &0.061 &0.014 &0.067 &0.355 &0.719 &0.812 &0.955
&0.440 &0.824 \\
 &$n=100$ &0.008 &0.043 &0.009 &0.046
&0.012 &0.050 &0.946 &0.987 &1.000 &1.000 &0.984 &0.997 \\
\hline
\end{tabular*}
\end{table}

Note that the condition that $(\Xvec-\muvec)/\|\Xvec-\muvec\|_2$ and
$(\muvec-\Xvec)/\|\muvec-\Xvec\|_2$ are identically distributed is
sufficient for the half-space median to have half-space depth 0.5, even
when $\Xvec$ lies in an arbitrary Banach space ${\mathcal B}$, where $\|
\cdot \|$ denotes the norm in ${\mathcal B}$. If $F$ is a probability
distribution over ${\mathcal B}$, and $\xvec$ is a fixed element in ${\mathcal
B}$, then the half-space depth of $\xvec$ can be defined as $\operatorname{HD}(\xvec,
F)= \inf_{h \in B^{*}} P\{h(\Xvec-\xvec) \ge 0\}$, where $h \dvtx {\mathcal B}
\rightarrow \mathbb{R}$ is a linear functional that belongs to the dual
space ${\mathcal B^{*}}$, $P$ stands for the probability measure on ${\mathcal
B}$ corresponding to $F$, and $\Xvec$ is a~random element in ${\mathcal B}$
having the distribution $F$. The point ${\muvec} \in {\mathcal B}$ is
called a \textit{half-space median} if $\operatorname{HD}(\muvec, F)=\sup_{\xvec \in
{\mathcal B}} \operatorname{HD}(\xvec, F)$. Instead of Banach spaces, if we work
with a Hilbert space ${\mathcal H}$, due to the Riesz representation
theorem and the reflexive nature of a Hilbert space, the half-space
depth of an observation $\xvec \in {\mathcal H}$ can be defined as
$\operatorname{HD}(\xvec,F) = \inf_{\hvec \in {\mathcal H}} P\{\langle\hvec,(\Xvec -\xvec)\rangle  \ge
0\}$, where $\langle\cdot,\cdot\rangle$ stands for the inner product defined on
${\mathcal H}$.

From the above discussion, it is clear that if we have a symmetric
distribution in a Hilbert or  Banach space, then the point of symmetry
will achieve the maximum depth value 0.5, and it will be the half-space
median. So, in a sense, the half-space median is well defined and
behaves in a~nice way, even in infinite-dimensional spaces for
symmetric probability distributions. However, in infinite-dimensional
spaces, even when we deal with nice symmetric distributions, the
half-space depth function can exhibit some anomalous behavior, which we
will see in the next section.

\section{Anomalous behavior of half-space depth in infinite-dimensional
spaces}\label{sec4}

We know that if we have a data cloud of $n$ observations in a
$d$-dimensional space, then the empirical depth of an observation lying
outside the convex hull formed by the data cloud is zero. For $d>n$,
since the Lebesgue measure of this convex hull is zero, we have zero
depth for all points in a set of probability measure one whenever we
have $n$ i.i.d. observations from an absolutely continuous distribution
in $\mathbb{R}^d$. In fact, for any probability measure on an
infinite-dimensional Banach space such that any finite-dimensional
hyperplane in that space has zero probability, the empirical half-space
depth based on finitely many i.i.d. observations from that probability
distribution will be zero almost everywhere. So, the empirical version
of half-space depth does not carry any statistically useful information
in such cases. Naturally, we would be curious to know what happens to
the population depth function in such situations. The following theorem
demonstrates that it is possible to have a nice symmetric probability
distribution on the~$l_2$ space for which the population depth function
takes positive values only on a set of probability measure zero. Recall
that the $l_2$ space of real sequences consists of infinite sequences
$(x_1,x_2,\dots)$  such that $\sum_{i=1}^{\infty} x_i^2 < \infty$.

\begin{theorem}\label{theo3}
Consider an infinite sequence of
independent random variables $\Xvec = (X_1, X_2, X_3, \ldots)$, where
$E( X_i)=0$ and $E(X_i^2) = \sigma_{i}^{2}$ for all $i\ge 1$ such that
$\sum_{i=1}^{\infty} \sigma_{i}^{2} < \infty$. Note that this implies
that $\Xvec$ lies in the $l_2$ space of real sequences with probability
one. Also, assume that the $X_i$'s have finite fourth moments and that
$\sum_{i=1}^{\infty} E(X_i^4)/i^2\sigma_i^4 < \infty$. For instance,
all these conditions will hold if the $X_i$'s are independent Gaussian
random variables. Then, for any given $\xvec=(x_1,x_2,\ldots)$ in that
$l_2$ space, the half-space depth of $\xvec$ with respect to the
distribution of $\Xvec$ will be zero unless $\xvec$ lies in a subset
having probability zero.
\end{theorem}

The proof of this theorem is given in the \hyperref[appendix]{Appendix}. This theorem
clearly shows that not only the empirical version, but also the
population version of the half-space depth will exhibit anomalous
behavior for some very common distributions in infinite dimensions.
Since any separable Hilbert space is isometrically isomorphic to the
$l_2$ space in view of the existence of a countable orthonormal basis
in such a space, similar examples can also be constructed on separable
Hilbert spaces. Clearly, the half-space depth function will not be a
very useful statistical concept in such spaces. To conclude, let us
recall the property of half-space depth characterizing the underlying
distribution established by earlier authors that was discussed in the
\hyperref[sec1]{Introduction}. From the above discussion, it is clear that in a
separable Hilbert space, there exist several probability measures,
which may even have independent Gaussian marginals, with half-space
depth functions identically equal to zero except on a subset having
zero probability measure. Nevertheless, such symmetric probability
measures will have a well-defined half-space median that achieves the
depth value 0.5.

\begin{appendix}\label{appendix}
\section*{Appendix}

\begin{lemma}\label{lem1}
Let $\operatorname{HD}(\xvec, F)$ be the half-space depth of
$\xvec$ with respect to the distribution $F$, and $F$ have density $f$
of the form $f(\xvec) = \psi(\|\xvec\|_p)$ with a monotonically
decreasing function $\psi$ and $0<p \le \infty$. Then, for any
$\xvec=(x,0,\ldots,0)$ on the coordinate axis, we have $\operatorname{HD}(\xvec,F) =
P(X_1 \geq x)$ when $x>0$, and $\operatorname{HD}(\xvec,F) = P(X_1 \leq x)$ when $x
\le 0$.
\end{lemma}

\begin{pf} We will prove it for
$\xvec_0=(1,0,\ldots,0)$. Proof for other points follows in the same
way. Consider any hyperplane $\alphavec(\xvec - \xvec_0)^{\prime}=0$ other
than $x_1=1$ that passes through $\xvec_0$ (see the left-hand diagram
in Figure \ref{fig1} for the case $d=2$). Here,
$\alphavec=(\alpha_1,\alpha_2,\ldots,\alpha_d)$ is a vector in
$\mathbb{R}^d$. Define the regions $A_1=\{\xvec\dvtx   x_1 < 1  \mbox{ and }
\alphavec(\xvec - \xvec_0)^{\prime} \ge 0\}$ and $A_2=\{\xvec\dvtx   x_1 \ge 1
  \mbox{ and }  \alphavec(\xvec - \xvec_0)^{\prime} < 0\}$ (see the left-hand
diagram in Figure \ref{fig1} for the case $d=2$). To prove the lemma, we have to
show that $P(\Xvec \in A_1) \ge P(\Xvec \in {A}_2)$. Define $A_3=\{
\xvec = (x_1,x_2,\ldots,x_d) \dvtx (x_1, -x_2, -x_3,\ldots, -x_d) \in
A_2\}$. Because of the symmetry of $f$, it is easy to check that
$P(\Xvec \in {A}_2)=P(\Xvec \in {A}_3)$. Therefore, it is enough to
prove that $P(\Xvec \in A_1) \ge P(\Xvec \in {A}_3)$. Note that for
every point $\zvec=(x_1,x_2,\ldots,x_d)$ in $A_1$, we have a point
$\zvec^{\prime}=(x_{1}^{\prime},x_2,x_3,\ldots,x_d)$ in $A_3$ such that
$x_1^{\prime}=2x_1-1$. Hence, $|x_1|\le |x_{1}^{\prime}|$ and $\|\zvec\|_p \le
\|\zvec^{\prime}\|_p$ with strict inequality being true for all $\zvec$ not
lying on the hyperplane $x_1=1$. This implies that $f(\zvec) \geq
f(\zvec^{\prime})$. Since the strict inequality holds over a set of positive
measure, integrating $f(\zvec)$ (resp. $f({\zvec}^{\prime})$) with respect
to $\zvec$ (resp. ${\zvec}^{\prime}$), we actually get  $P(\Xvec \in A_1)> 
P(\Xvec \in {A}_3)$.\looseness=-1
\end{pf}

\begin{lemma}\label{lem2}
Consider a p.d.f. $f$ on
$\mathbb{R}^d$ satisfying $f(\xvec) = \psi(\|\xvec\|_\infty)$ and a
random vector $\Xvec$ with  p.d.f. f. Then, for any $x>0$, we have
$P(X_1+X_2 \ge 2x) < P(X_1 \ge x)$.
\end{lemma}

\begin{pf} Again, we will prove this only
for $x=1$. Let us define $A_1=\{\xvec=(x_1,x_2,\ldots,x_d) : x_1<1
  \mbox{ and }  x_1+x_2\ge 2\}$ and $A_2= \{ \xvec=(x_1,x_2,\ldots,x_d) :
x_1\ge 1   \mbox{ and }  x_1+x_2<2\}$ (these two regions are shown in the
right-hand diagram in Figure \ref{fig1} for the case $d=2$). We also define the
region $A_3=\{\xvec=(x_1,x_2,\ldots,x_d) \dvtx  (x_2,x_1,x_3,\ldots,x_d)
\in A_1\}$. Because of the symmetry of $f(\xvec)$ under permutations of
the coordinates of $\xvec$, it is straightforward to see that $P(\Xvec
\in A_1)=P(\Xvec \in A_3)$. Hence, it is enough to show that $P(\Xvec
\in A_3) < P(\Xvec \in A_2)$. Now, for any $\zvec=(z_1,z_2,\ldots,z_d)
\in A_2$, we have a corresponding point $\zvec^{\prime}=(2-z_2,2-z_1,z_3,
\ldots,z_d)$ in $A_3$. Also, note that for any
$\zvec=(z_1,z_2,\ldots,z_d)$ in $A_2$, $z_1$ and $z_2$ have the
respective forms $z_1=1+b$ and $z_2=1-b-a$ for some $a,b >0$ (see the
right-hand diagram in Figure \ref{fig1} for the case $d=2$). Consequently, for
$\zvec^{\prime}=(z_{1}^{\prime},z_{2}^{\prime},z_3,\ldots,z_d)$, we have
$z_{1}^{\prime}=1+b+a$ and $z_{2}^{\prime}=1-b$. Clearly, $\max\{|z_1|,|z_2|\} <
\max\{ |z_{1}^{\prime}|,|z_{2}^{\prime}|\}=1+a+b$, which implies that
$\|\zvec\|_\infty \leq \|\zvec^{\prime}\|_\infty$ and hence that $f(\zvec) >
f(\zvec^{\prime})$ with strict inequality on a set of positive probability
measure under $f$. This proves that $P(\Xvec \in A_2)>P(\Xvec \in
A_3)$.
\end{pf}

\begin{lemma}\label{lem3}
 Let
$f(\xvec)=\psi(\|\xvec\|_p)$ for $1\le p<\infty$ be the p.d.f. of
$\Xvec=(X_1, X_2, \ldots, $ $X_d)$. Consider $\xvec_0=(c,c,0,\ldots,0)$
for $c>0$. Its half-space depth is then given by $\operatorname{HD}(\xvec_0,F)$ =
$P(X_1+X_2 \geq 2c)$.
\end{lemma}

\begin{pf}  Consider the hyperplane
$x_1+x_2=2c$ (see Figure \ref{fig2} for the case $d=2$). We have to show that
this hyperplane determines the half-space depth of $\xvec_0$. For this,
we will follow the same lines of argument as in Lemmas \ref{lem1} and \ref{lem2}.
Consider a new hyperplane $\alphavec(\xvec-\xvec_0)^{\prime}=0$ passing
through $\xvec_0$ (see Figure \ref{fig2} for the case $d=2$). Define the regions
$A_1=\{\xvec=(x_1,x_2,\ldots,x_d)\dvtx   x_1+x_2 < 2c  \mbox{ and }
\alphavec(\xvec - \xvec_0)^{\prime} \ge 0\}$ and
$A_2=\{\xvec=(x_1,x_2,\ldots,x_d)\dvtx   x_1 +x_2\ge 2c   \mbox{ and }
\alphavec(\xvec - \xvec_0)^{\prime} < 0\}$ (see Figure \ref{fig2} for the case
$d=2$). To prove the lemma, we have to show that $P(\Xvec \in A_1) \ge
P(\Xvec \in {A}_2)$. Define $A_3=\{ \xvec = (x_1,x_2,\ldots,x_d) \dvtx
(x_2, x_1, x_3,\ldots, x_d) \in A_2\}$. Because of the symmetry of
$f(\xvec)$ under any permutation of the coordinates of $\xvec$, we have
$P(\Xvec \in A_2)=P(\Xvec \in A_3)$. Therefore, it is enough to show
that $P(\Xvec \in A_3) \le P(\Xvec \in A_1)$.

Note that any point $\zvec \in A_1$ is of the form
$\zvec=(c+a,c-a-k,x_3,\ldots,x_d)$, where $k>0$, and $a$ can be
positive or negative (see Figure \ref{fig2} for the case $d=2$). For any $\zvec
\in A_1$, we get a~corresponding point $\zvec^{\prime} \in A_3$ such that
$\zvec^{\prime}=(c+a+k, c-a,x_3,\ldots,x_d)$. We now need to show that
$\|\zvec\|_p \le \|\zvec^{\prime}\|_p$ and for that, we will consider the
two cases $a>0$ and $a<0$ separately.

When $a>0$ (see the left-hand diagram in Figure \ref{fig2}
 for the case $d=2$),
we have $0 < |c-a| < |c+a|$. Now, for $p\ge 1$ and $t,k>0$, it is easy
to check that the function $h(t) = (t+k)^p - t^p$ is non-decreasing in
$t$. So, for $0<t_1<t_2$, we have $0 < h(t_1)\le h(t_2)$. Taking
$t_1=|c-a|$ and $t_2=|c+a|$, we get $(|c-a|+k)^p - |c-a|^p \le
(|c+a|+k)^p - |c+a|^p$. Now, using the facts that $ |c+a|+k=|c+a+k|$
and $|c-a-k| \leq |c-a|+k$, we arrive at $|c-a-k|^p - |c-a|^p \le
|c+a+k|^p - |c+a|^p$. This implies that $|c-a-k|^p + |c+a|^p \le
|c+a+k|^p + |c-a|^p$, which in turn implies that $\|\zvec\|_p \le
\|\zvec^{\prime}\|_p$. Note that the strict inequality holds on a set of
positive probability measure under $f$.

For $a<0$ (see the right-hand diagram in Figure \ref{fig2} in the case $d=2$),
first note that $a+k>0$ and that the coordinates of $\zvec$ and
$\zvec^{\prime}$ are of the respective forms $\zvec=(c-\alpha,
c-\beta,x_3,\ldots,x_d)$ and $\zvec^{\prime}=(c+\alpha, c+
\beta,x_3,\ldots,x_d)$, where $\alpha =-a>0$ and $\beta=a+k>0$. Now,
$|c-\alpha|<|c+\alpha|$ and $|c-\beta|<|c+\beta|$ imply that
$\|\zvec\|_p < \|\zvec^{\prime}\|_p$.
\end{pf}

\begin{lemma}\label{lem4}
Assume that we have a
p.d.f. $f$ that satisfies $f(\xvec)=\psi(\|\xvec\|_p)$ for some $p>0$
and monotonically decreasing $\psi$. Let $\Xvec=(X_1,X_2,\ldots,X_d)$
be a random vector with p.d.f. f. If $X_1$ and $2^{(1-p)/p}(X_1+X_2)$
are identically distributed, then we must have $p=2$.
\end{lemma}

\begin{pf} First, note that if
$f(\xvec)=\psi(\|\xvec\|_p)$, then the joint p.d.f. of $X_1$ and $X_2$
is of the form $f_1(x_1,x_2)=\psi_1(|x_1|^p+|x_2|^p)$ for some $\psi_1
\dvtx \mathbb{R}_+ \rightarrow \mathbb{R}_+$. We can show that the p.d.f.'s
of $X_1$ and $Y=2^{\alpha}(X_1+X_2)$, where $\alpha=(1-p)/p$, are given
by $f_{X_1}(x)=\int \psi_1((|x|^p+|x_2|^p)^{1/p})\,\mathrm{d}x_2$ and
$f_{Y}(x)=2^{-\alpha}\int \psi_1((|2^{-\alpha}x-x_2|^p+|x_2|^p)^{1/p})\,\mathrm{d}x_2$, respectively.

Since both of these  p.d.f.'s are continuous functions, and $X_1$ and
$Y$ are identically distributed, we can equate their values at $x=0$.
We then get $\int \psi_1(|x_2|)\,\mathrm{d}x_2 = 2^{-\alpha} \int
\psi_1(2^{1/p}|x_2|)\,\mathrm{d}x_2= 2^{-(\alpha+1/p)}\int \psi_1(|x_2|)\,\mathrm{d}x_2$.
Hence, we must have $\alpha=-1/p$, which implies $p=2$.
\end{pf}

\begin{pf*}{Proof of Theorem \ref{theo2}} Note that the ``if''
part is trivial in view of our discussion preceding the statement of
the theorem. We shall now prove the ``only if'' part.

First, we shall prove it for the bivariate case, that is, $d=2$.
Without loss of generality, we assume that $\muvec={\mathbf 0}$. Let $Z$ be
the angle between the positive side of the $x_1$-axis and the random
vector $\Xvec$ (measured counterclockwise from the $x_1$-axis). Now,
consider a straight line which passes through the origin and makes an
angle $\theta$ with the $x_1$-axis. Since $\muvec={\mathbf 0}$, the two
half-spaces generated by that straight line will have the same
probability measure. Now, rotate the line in a~counterclockwise
direction by an angle $\delta$ to bring it to a new position. Clearly,
the two half-spaces generated by the straight line in the new position
will also have the same probability~0.5. This implies that $P(\theta <
Z < \theta+\delta) = P(\uppi+ \theta < Z < \uppi+ \theta +\delta)$. Since
this equality holds for all $\theta$ and $\delta$, it implies that $Z$
and $Z+\uppi$ have the same probability distribution. The result now
follows from the fact that $(\Xvec -\muvec)/\|\Xvec - \muvec\|_2 =(\operatorname{Cos}
Z,   \operatorname{Sin} Z)$ and $(\muvec -\Xvec)/\|\Xvec - \muvec\|_2 =(\operatorname{Cos} (Z+\pi),
  \operatorname{Sin} (Z+\pi) )$.

For $d>2$, we need to consider $d-1$ random angles
$Z_1,Z_2,\ldots,Z_{d-1}$. Note that here the direction vector $(\Xvec
-\muvec)/\|\Xvec - \muvec\|_2$ can be expressed as $(\Xvec
-\muvec)/\|\Xvec - \muvec\|_2$= $(\operatorname{Cos} Z_1,\break \operatorname{Sin} Z_1 \operatorname{Cos}Z_2,  \ldots ,
\operatorname{Sin} Z_1\cdots \operatorname{Sin} Z_{d-2} \operatorname{Cos} Z_{d-1},  \operatorname{Sin} Z_1\cdots \operatorname{Sin}
Z_{d-2} \operatorname{Sin}Z_{d-1})$. Now, consider a~hyperplane $H$ which makes angles
$\theta_1,\theta_2,\ldots,\theta_{d-1}$ with the coordinate axes and
then rotate it to $H_1$ such that the new angles are
$\theta_1+\delta,\theta_2,\ldots,\theta_{d-1}$. The result now follows
from the same argument that is used in the bivariate case.
\end{pf*}

\begin{lemma}\label{lem5}
For any two sequences
${\sigmavec}=(\sigma_1,\ldots)$ and $\xvec=(x_1, x_2, \ldots)$ in the
$l_2$ space of real sequences, we have $\sup_{\alphavec \in l_2}
 \{ ({\sum_{i=1}^{\infty}{\alpha_{i}^2 \sigma_i^{2}}})^{-1/2}
({\sum_{i=1}^{\infty} \alpha_i x_i} )  \} < \infty$ if and only if
$\sum_{i=1}^{\infty} x_i^2/\sigma_i^2 < \infty$.
\end{lemma}

\begin{pf} (The ``if'' part). For any
$\alphavec \in l_2$, ${\sum_{i=1}^{\infty} \alpha_i x_i} \le
({\sum_{i=1}^{\infty} \alpha_{i}^2\sigma_{i}^{2}})^{1/2}
({\sum_{i=1}^{\infty} x_i^2/\sigma_i^2})^{1/2}$ (i.e., the
Cauchy--Schwarz inequality) implies that ${\sum_{i=1}^{\infty} \alpha_i
x_i}/ ({\sum_{i=1}^{\infty}{\alpha_{i}^2 \sigma_i^{2}}})^{1/2} \le
\sum_{i=1}^{\infty} x_i^2/\sigma_i^2$. Now, the right-hand side of the
inequality does not depend on $\alphavec$. So, $\sum_{i=1}^{\infty}
x_i^2/\sigma_i^2<\infty$  implies the finiteness of $\sup_{\alphavec
\in l_2} \{ {\sum_{i=1}^{\infty} \alpha_i x_i}/
({\sum_{i=1}^{\infty}{\alpha_{i}^2 \sigma_i^{2}}})^{1/2} \} \le
\sum_{i=1}^{\infty} x_i^2/\sigma_i^2$.

(The ``only if'' part). Next, consider the case where
$\sum_{i=1}^{\infty} x_i^2 / \sigma_i^2 = \infty$. Choose a sequence
$\{\alphavec_n\}$ of real sequences, where
$\alphavec_n=(\alpha_{n1},\alpha_{n2},\ldots)$ has non-zero values only
at first $n$ coordinates (i.e., $\alpha_{ni}=0$ for all $i>n$) and
$\alpha_{ni}=x_i/\sigma_i^2$ for $i=1,2,\ldots,n$. Clearly,
$\alphavec_n \in l_2$ for all $n \geq 1$, and for each $n$, it is easy
to check that ${\sum_{i=1}^{n} \alpha_{ni}  x_i}/
({\sum_{i=1}^{n}{\alpha_{ni}^2 \sigma_i^{2}}})^{1/2}=({\sum_{i=1}^{n}
x_i^2/\sigma_i^2})^{1/2}$. So, we get $\sup_{n \geq 1}
 \{{\sum_{i=1}^{n} \alpha_{ni}  x_i}/
({\sum_{i=1}^{n}{\alpha_{ni}^2 \sigma_i^{2}}})^{1/2} \} =\infty$.
This clearly implies that we have $\sup_{\alphavec \in l_2} \{
{\sum_{i=1}^{\infty} \alpha_{i}x_i}/({\sum_{i=1}^{\infty}{\alpha_{i}^2 \sigma_i^{2}}})^{1/2}\} $ $
=\infty$.
\end{pf}

\begin{pf*}{Proof of Theorem \ref{theo3}} Consider any $\xvec$
in the $l_2$ space with $\xvec \neq {\mathbf 0}$. For any $\alphavec$ in
the $l_2$ space, the random variable $Z=\langle\alphavec,\Xvec\rangle$ has a
probability distribution with $E(Z)=0$ and $V(Z)=\sum_{i=1}^{\infty}
\alpha_{i}^2 \sigma_i^{2}$. Using Chebyshev's inequality, we get
$P(\langle\alphavec,(\Xvec-\xvec)\rangle  \ge 0)=P(Z  \ge  \langle\alphavec,\xvec\rangle) \le
\sum_{i=1}^{\infty} \alpha_{i}^2 \sigma_i^{2}/(\sum_{i=1}^{\infty}
\alpha_i x_i)^2.$ So, the depth of $\xvec$ is bounded above by
$\inf_{\alphavec \in l_2} \{\sum_{i=1}^{\infty} \alpha_{i}^2
\sigma_i^{2}/(\sum_{i=1}^{\infty} \alpha_i x_i)^2 \}$. From Lemma
\ref{lem5}, it follows that this upper bound is zero when $\sum_{i=1}^{\infty}
x_i^2/\sigma_i^2 = \infty$. Therefore, $\xvec$ will have  positive
depth only if $\sum_{i=1}^{\infty} x_i^2/\sigma_i^2 \,{<}\,
\infty$.\looseness=-1

Next, consider $Y_i=X_i^2/\sigma_i^2$ for $i \ge 1$. The $Y_i$'s are
then independent random variables with a common mean $1$ and
$\sum_{i=1}^{\infty} E(Y_i^2)/i^2 < \infty$. So, using the strong law
of large numbers (see Theorem 1 in \cite{Chow}, page 124),
we have $n^{-1}\sum_{i=1}^{n} Y_i \stackrel{\mathrm{a.s.}}{\longrightarrow} 1$
as $n \rightarrow \infty$. Consequently, $\sum_{i=1}^{\infty} Y_i =
\sum_{i=1}^{\infty} X_i/\sigma_{i}^{2} = \infty$ with probability one.
\end{pf*}
\end{appendix}

\section*{Acknowledgements}

The authors are thankful to two anonymous referees for their careful
reading of an earlier version of the paper and for providing them with
several helpful comments. The first author would also like to thank
Professor B.V. Rao for helpful discussions.

\printhistory

\end{document}